\documentclass[12pt]{article}
\usepackage{amsmath,amsfonts,amscd,amssymb,latexsym}

\title{Algebraic Geometry over Free Metabelian  Lie Algebra II: Finite Field Case}
\author{\textsf{E. Yu. Daniyarova}\thanks{The first author is supported by
the RFFI grant N02-01-00192.} \and \textsf{I. V.
Kazachkov}\thanks{The second and the third authors are supported
by the `Universitety Rossii' grant.} \and \textsf{V. N.
Remeslennikov}$^{\dag}$}

\newtheorem{lem}{Lemma}[section]
\newtheorem{thm}{Theorem}[section]
\newtheorem{cor}{Corollary}[section]
\newtheorem{prop}{Proposition}[section]

\newtheorem{defn}{Definition}[section]
\newtheorem{expl}{Example}[section]
\newtheorem{rem}{Remark}[section]

\newcommand{\M}{\mathfrak{M}}
\newcommand{\V}{\mathfrak{V}}
\newcommand{\N}{\mathfrak{N}}
\newcommand{\F}{F}

\newcommand{\proof}{\paragraph{Proof.}}
\renewcommand{\hom}{\texttt{Hom}}
\newcommand{\rad}{\texttt{Rad}}
\newcommand{\id}{\texttt{id}}
\newcommand{\AS}{\texttt{AS}}

\newcommand{\AX}{A \left[ X \right]}

\newcommand{\sFit}{\textsf{Fit}}
\newcommand{\Fit}{\texttt{Fit}}
\newcommand{\ucl}{\texttt{ucl}}

\begin{document}

\maketitle

\begin{abstract}
This paper  is the second in a series of three, the aim of
which is to construct algebraic geometry over a free  metabelian Lie
algebra $F$. For the universal  closure of free metabelian Lie
algebra of finite rank $r \ge 2$ over a finite field $k$ we find a
convenient set of axioms in the language of Lie algebras $L$ and the language $L_{F}$ enriched by constants from $F$. We give a description of:
\begin{enumerate}
    \item the structure of finitely generated algebras from the universal
        closure of $F_r$ in both $L$ and $L_{F_r}$;
    \item the structure of irreducible algebraic sets over $F_r $ and respective
        coordinate algebras.
\end{enumerate}
We also prove that the universal theory of a free metabelian Lie algebra over a finite field is decidable in both languages.
\end{abstract}

\tableofcontents

\section{Introduction} \label{sec:1}
This paper is the second in a series of papers the main object of
which is to construct algebraic geometry over free metabelian Lie
algebra. In this paper we consider the free metabelian Lie algebra
$F_r $ of a finite rank $r \ge 2$ over a finite field $k$.
Throughout this paper we use the results, notation and definitions
of the first paper of the current series \cite{AGFMLA1}.

The object of Section \ref{sec:2}, which arises from papers
\cite{alg1} and \cite{alg2} is to lay the foundations of algebraic
geometry over Lie algebras. In \cite{alg1} and \cite{alg2} the
authors conduct their arguments and prove the results in the
category of groups, however the proofs are absolutely analogous
for Lie algebras and, therefore, most of the results in Section
\ref{sec:2} are omitted.

Section \ref{sec:3} holds main technical complications of the
current paper. There we introduce two collections of seven series
of universal axioms $\Phi _r$ and $\Phi'_r $, $r \ge 2$. The
axioms of the collection $\Phi _r$ are universal formulas of the
standard first order language $L$ of theory of Lie algebras over
the field $k$ and the axioms of the collection $\Phi'_r $ are
universal formulas in the enriched language $L_{F_r}$, obtained
from $L$ by joining constants for the elements of $F_r$. There we
establish some properties of Lie algebras that satisfy either of
the collections mentioned.

The key results of the current paper are formulated in Section
\ref{sec:4} (see Theorems \ref{thm:41} - \ref{thm:45}). We list
the most important of these results:
\begin{itemize}
    \item for the universal closure of the free
    metabelian Lie algebra of finite rank $r \ge 2$ over a finite
    field $k$ we find a convenient set of axioms ($\Phi_r $ and $\Phi '_r $) in $L$ and $L_{\F}$,
    \item describe the structure of finitely
          generated algebras from  $\F_r- \ucl (\F_r) $ and $\ucl
          (\F_r)$,
    \item prove that the universal theory of the free metabelian Lie algebra over a finite field is decidable
    in both $L$ and $L_{\F}$.
\end{itemize}

In Section \ref{sec:5} we apply theorems from Section \ref{sec:4}
to algebraic geometry over the algebra $F_r $, $r \ge 2$ over a
finite field $k$. The main results of this section are:
\begin{itemize}
    \item given a structural description of coordinate algebras of irreducible
    algebraic sets over  $F_r $;
    \item given a description of the structure of irreducible
    algebraic sets;
    \item constructed a theory of dimension in the category of
    algebraic sets over  $F_r $.
 \end{itemize}

\section{Elements of Algebraic Geometry over Lie Algebras}
\label{sec:2}

In paper \cite{alg1} the authors introduce main notions of
algebraic geometry over groups. In Subsection \ref{sec:21} below
we introduce main notions of algebraic geometry over Lie algebras.
Following paper \cite{alg1} we list several results and theorems,
involving these notions. Subsection \ref{sec:22} highlightens some
of the aspects of algebraic geometry over the free metabelian Lie
algebra
 $F_r $ of finite rank $r$, $r \ge 2$

\subsection{General Case} \label{sec:21}

Let $A$ be a fixed Lie algebra over a field $k$.

Recall that a Lie algebra $B$ over a field $k$ is called an
$A$--Lie algebra if and only if it contains  a designated copy of
$A$, which we shall for most part identify with $A$. A
homomorphism $\varphi$ from an $A$--Lie algebra $B_1$ to an
$A$--Lie algebra $B_2$ is an $A$-homomorphisms of Lie algebras if
it is the identity on $A$, $\varphi (a)=a, \; \forall a \in A$.
Set $\hom _A (B_1, B_2)$ to be the set of all $A$--homomorphisms
from $B_1$ to $B_2 $. We use the symbol `$\cong _A$'
($A$--isomorphism) to express that two $A$--Lie algebras are
isomorphic in the category of $A$--Lie algebras.

 The family of all $A$--Lie
algebras together with the collection of all $A$--homomorphisms
form a category in the obvious way.

Let $X = \left\{ x_1, \ldots, x_n \right\} $ be a finite set. The
free $A$--Lie algebra with the free base $X$
$$
\AX = A \ast F(X),
$$
is the free Lie product of the free (in the category of Lie
$k$-algebras) Lie algebra $F(X)$ and the algebra $A$. We think of
elements of $\AX$ as polynomials with coefficients in $A$. We use
functional notation here,
$$
f = f(x_1 ,\ldots ,x_n ) = f(x_1 ,\ldots ,x_n ,a_1 ,\ldots ,a_r )
$$
thereby expressing the fact that the Lie polynomial representing
$f$ in $\AX$ involves the  variables $x_1,\dots,x_n$ and, as
needed, the constants $a_1,\dots,a_r \in A$.

Using the standard argument from universal algebra one verifies
that the algebra $\AX$ is the free algebra in the category of
$A$--algebras.

Let $B$ be an $A$--Lie algebra and let $S$ be a subset of $\AX$.
Then the set
$$
B^n =  \left\{ (b_1 ,\ldots ,b_n ) | \;b_i  \in B \right\} .
$$
 is termed the affine $n$-dimensional space over the algebra $B$.

A point  $p = (b_1,\ldots , b_n ) \in B^n $ such that
$$
 f(p) = f(b_1 ,\ldots ,b_n, a_1 ,\ldots,a_r ) = 0.
$$
is termed a  root of the polynomial $f \in \AX$. In that case we
also say that the polynomial $f$ vanishes at the point $p$.

A point $p \in B^n $ is  a root or a solution of the system $S
\subseteq \AX$  if every polynomial from $S$ vanishes at $p$, i.e.
if $p$ is a root of every polynomial from the system $S$.

The set
 $$
V_B (S) = \left\{ p \in B^n  | \; f(p) = 0 \ \ \forall \, f \in S
\right\}
$$
is termed the (affine) algebraic set over $B$ defined by the
system of equations $S$.

Let $S_1$ and $S_2$ be subsets of $\AX$. Then the systems $S_1$
and $S_2$ are termed equivalent over $B$ if \/ $V_B (S_1 ) = V_B
(S_2)$.

\begin{expl}[typical examples of algebraic sets] \label{expl:1} \
\begin{enumerate}
    \item Every element $a \in A$ forms an algebraic set,
 $\left\{ a \right\}$: $S = \left\{ x - a = 0 \right\}$,
 $V_B (S) = \left\{ a \right\}$. In this example $n = 1$
 and $X = \left\{ x \right\}$.
       \item The centraliser $C_B (M)$ of an arbitrary set of elements $M$ from
$A$ is the algebraic set defined by the system $S = \left\{ x
\circ m = 0|\; m \in M \right\}$.
\end{enumerate}
\end{expl}

Let $Y$ be an arbitrary   algebraic set ($Y = V_B (S)$) from $B^n
$. The set {\rm
$$
\rad_B (S)=\rad_B (Y) =  \left\{ f \in \AX |\; f(p) = 0 \ \
\forall \, p \in Y \right\}
$$}
is termed the radical of the set $Y$. If \/ $Y=\emptyset$ then, by
the definition, its radical is the algebra $\AX$.

Clearly, the radical of a set is an ideal of the algebra $\AX$. A
polynomial $f \in \AX$ is termed a consequence of a system $S
\subseteq \AX$ if
$$
V(f) \supseteq V(S).
$$
The radical of an algebraic set describes it uniquely, i.e. for
two arbitrary algebraic sets $Y_1 ,Y_2  \subseteq B^n $ {\rm
$$
Y_1  = Y_2 \hbox{ if and only if }  \rad_B (Y_1 ) = \rad_B (Y_2 ).
$$}

Let $B$ be an $A$--Lie algebra, $S$ be a subset of $\AX$ and
$Y\subseteq B^{n}$ be the algebraic set defined by the system $S$.
Then the factor-algebra {\rm
$$
\Gamma _B (Y) = \Gamma _B (S) = {\raise0.7ex\hbox{${\AX}$}
\!\mathord{\left/
 {\vphantom {{\AX} {\rad_B (Y)}}}\right.\kern-\nulldelimiterspace}
\!\lower0.7ex\hbox{${\rad_B (Y)}$}}
$$}
is termed the coordinate algebra of the algebraic set $Y$ (or of
the system $S$).

Note that coordinate algebras of  consistent  systems of equations
are $A$--Lie algebras and form a subcategory of the category of
all $A$--Lie algebras.

\begin{lem} \label{prop:algsashom} \label{lem:211}
For any algebraic set $V_B(S)$ there is a one-to-one
correspondence between the points of $V_B(S)$  and
$A$--homomorphisms from $\Gamma _B(V_B(S))$ to $B$.
\end{lem}
\proof To a point $y \in V_B(S)$ we link an $A$--homomorphism
$$
\varphi \in \hom_A (\AX,B) \hbox{ given by } f \rightarrow f(y), \
f \in \AX.
$$
This homomorphism on $\AX$ is a correct homomorphism on the
factor-algebra $\Gamma _B (S)$. Conversely, if $\varphi \in \hom_A
(\AX,B)$ then the point corresponding to $\varphi$ is the
following $y=(\varphi(\bar x_1),\ldots, \varphi(\bar x_n))$, where
$\bar x_1,\ldots, \bar x_n$ are the images of  $X=\left\{x_1
,\ldots ,x_n\right\}$ in the factor-algebra $\Gamma _B(V_B(S))$.
Clearly, $y \in V_B(S)$. Obviously, the two maps given are
mutually inverse.\hfill $\blacksquare$

\begin{rem}
The coordinate algebra  $\Gamma (Y)$ is an $A$--Lie algebra
isomorphic to the algebra of all polynomial functions from  $Y$
into $B$ defined by the rule $y \in Y$, $y \rightarrow f(y)$ for
$f \in \AX$.
\end{rem}

\begin{expl}[the coordinate algebra of a point]
\label{ex:coordalg}
If $a \in A$, $Y = \left\{ a \right\}$ then $\Gamma _A (Y) \cong
A$.
\end{expl}

We next fix an $A$--Lie algebra $B$. Let $Y \subseteq B^n$ and $Z
\subseteq B^m$ be algebraic sets. Then the map
$$
\psi: Y \longrightarrow Z
$$
is termed a morphism from the algebraic set $Y$ to the algebraic
set $Z$ if there exist $f_1,\dots,f_m \in A\left[x_1,\ldots ,x_n
\right]$ so that for any $(b_1,\ldots,b_n) \in Y$
$$
\psi(b_1,\dots,b_n)=(f_1(b_1,\ldots,b_n),\ldots,f_m(b_1,\ldots,b_n))
\in Z.
$$
Algebraic sets $Y$ and $Z$  are termed isomorphic if there exist
morphisms:
$$\psi: Y \rightarrow Z \hbox{ and } \theta: Z \rightarrow Y$$
such that $\theta \psi  = \verb"id"_Y $ and \/ $\psi \theta  =
\verb"id"_Z$. We shall make use of the notation $\hom(Y,Z)$ for
the set of all morphisms from $Y$ to $Z$.

The collection of all algebraic sets $V_B (S)$, $S \subset
A\left[x_1 ,\ldots,x_n \right]$ over $B^{n}$, where  $n$ is a
non-fixed positive integer form the family of objects of
$\AS_{A,B}$. Morphisms of this category are the morphisms of
algebraic sets.

Following the argument of paper \cite{alg1} one can prove that the
categories of coordinate algebras and algebraic sets are
equivalent. We formulate this result by the means of the following
two lemmas:

\begin{lem} \label{lem:212}
Coordinate algebras define algebraic sets up to isomorphism:
\[
Y \cong Y'\quad  \Leftrightarrow \quad \Gamma (Y) \cong _A \Gamma
(Y').
\]
\end{lem}

\begin{lem} \label{lem:213}
There exists a one-to-one correspondence between {\rm
$\hom(Y,Y')$} and {\rm $\hom_A (\Gamma (Y'),\Gamma (Y))$}.
Furthermore, whenever we have an embedding of algebraic sets $ Y
\subseteq Y'$ the correspondent map $\varphi :\Gamma (Y')
\rightarrow \Gamma (Y)$ is an $A$--epimorphism  of coordinate
algebras. Moreover, if $Y\subsetneq Y_1$ then the kernel
$\ker\;\varphi \ne 0$ is non-trivial.
\end{lem}

\begin{expl} \label{expl:3} \label{expl:I-dim1}
A vector space $A$ over a field $k$ with the trivial
multiplication {\rm(}$\forall u,v \in A \, \left[u\circ
v=0\right]${\rm)} is a particular case of a Lie $k$-algebra.
 Let $A$ be a Lie $k$-algebra with the trivial multiplication and
assume that $B = A$. Applying theorems of linear algebra one shows
that:
\begin{enumerate}
    \item Every consistent system of equations over $A$ is equivalent
to a triangular system of equations (see \cite{LinAlg} for
definitions).
    \item Every algebraic set $Y \subseteq A^n$ is
isomorphic to an algebraic set of the form
$$
(\underbrace {A,A,\ldots,A}_s,0,\ldots,0), \ \ 0 \le s \le n.
$$
    \item  Every coordinate algebra $\Gamma (Y)$ is  $A$--isomorphic to
$A \oplus \verb"lin"_k \left\{ x_1 ,\ldots,x_s \right\}$, where
here $0 \le s \le n$, and \/ $\verb"lin"_k \left\{ x_1 ,\ldots,x_s
\right\}$ is the linear span of the elements $\left\{ x_1
,\ldots,x_s \right\}$ over $k$.
\end{enumerate}
\end{expl}

The union of two algebraic set is not necessarily again an
algebraic set. The correspondent counterexample is easy to
construct using Example \ref{expl:3}. We, therefore, define a
topology in $B^n$ by taking algebraic sets in $B^n$ as a sub-basis
for closed sets.  We term this topology the Zariski topology.

A closed set $Y$ is termed irreducible if \/ $Y=Y_1 \cup Y_2$,
where $Y_1$ and $Y_2$ are closed, implies that either $Y=Y_1$ or
$Y=Y_2$.

An $A$--Lie algebra $B$ is termed $A$--equationally Noetherian if
for every $n \in \mathbb{N}$ and for every system
 $S \subseteq A\left[x_1 ,\ldots,x_n \right]$
 there exists a finite subsystem  $S_0  \subseteq S$
such that $V_B (S) = V_B (S_0 )$.

\begin{thm} \label{thm:214}
Every closed subset $Y$ of $B^n$ over $A$--equationally Noetherian
$A$--Lie algebra $B$ can be expressed as a finite union of
irreducible algebraic sets:
$$
Y = Y_1 \cup \ldots \cup Y_l.
$$
These sets are usually referred to as the irreducible components
of $Y$, which turn out to be unique, if for every $i,j=1, \ldots
,l;\; i\ne j$ $Y_i \nsubseteq Y_j$.
\end{thm}

The main aim of algebraic geometry over an $A$--Lie algebra $B$ is
to give a description of algebraic sets over $B$ up to
isomorphism. Or, which is equivalent, to give a description of
coordinate algebras of algebraic sets up to $A$--isomorphism.

We next treat some of the properties of model-theoretical classes
generated by an $A$--Lie algebra $B$. Our interest to the
universal closure generated by $B$ is justified by the following
circumstance: finitely generated Lie algebras from the universal
closure $A - \ucl(B)$ are exactly the coordinate algebras of
irreducible algebraic sets over $B$.

Recall that the universal closure $A - \ucl(B)$ generated by $B$
is the class of all $A$--Lie algebras that satisfy all the
universal sentences satisfied by $B$ (for details see
\cite{AGFMLA1}).

\begin{thm} \label{thm:I-coor=ucl} \label{thm:215} \label{thm:222}
Let $B$ be an $A$--equationally Noetherian $A$--Lie algebra. Then
finitely generated $A$--Lie algebra $C$ is the coordinate algebra
of an irreducible algebraic set over $B$ if and only if {\rm $C
\in A - \ucl(B)$}.
\end{thm}

\subsection{The Case of $F_r$} \label{sec:22}

Let $\F_r$ be the free metabelian Lie $k$-algebra of the rank $r$,
let $\left\{ a_1,\ldots, a_r \right\}$ be its free base and let $R
= k\left[x_1 ,\ldots,x_r \right]$ be the ring of polynomials from
$r$ variables. Recall that the Fitting's radical of the algebra
$F_r$ coincides with its commutant $F_r^2$ and admits the
structure of an $R$--module. Further, the multiplication by the
variables $x_i$'s of the ring $R$ is interpreted as the
multiplication by free generators $a_i$'s (see \cite{AGFMLA1}).

In the current paper we consider so called `diophantine geometry',
i.e. we consider systems of equations with coefficients in  $\F_r$
and solutions of these systems from $\F_r$. In the event that
$r=1$, $\F_r$ is Abelian and this extreme case has been already
considered in Example \ref{expl:I-dim1} and we, therefore,
consider only non-degenerated alternative of $r \ge 2$.

One of the most important algebraic sets over the free metabelian
Lie algebra $\F_r$ is $\Fit(\F_r )$ is an algebraic set in the
affine space $\F_r ^1$. To prove this consider an equation $(a_1
a_2 )x = 0$ with one indeterminant $x$, where here $a_1$ and $a_2
$ are two distinct elements of the free base of $\F_r $. Since
$\F_r $ is a $U$-algebra (see \cite{AGFMLA1}), its Fitting's
radical is Abelian. Consequently every element of $\Fit(\F_r )$
satisfies this equation and $\Fit(\F_r ) \subseteq
V_{\F_r}(\left\{(a_1 a_2 )x = 0\right\})$. To prove the reverse
inclusion  take $c \notin \Fit(\F_r )$. Since $\Fit(\F_r )$ is a
torsion free module over the ring of polynomials $R$, we obtain
$(a_1 a_2 )c \ne 0$, which implies that $\Fit(\F_r ) =
V_{\F_r}(\left\{(a_1 a_2 )x = 0\right\})$. Below we show that the
Fitting's radical is an irreducible algebraic set and that its
coordinate algebra is $\F_r \oplus T_1 $, where here $T_1$ is the
$R$-free module of the rank 1 (for definition of $\F_r \oplus
T_s$, where here $T_s$ is the $R$-free module of the rank $s \ge
1$, see \cite{AGFMLA1}).
\medskip
\newline
\textbf{Notation.} {\it By $\F_{r,s}$ we denote the direct module
extension of the Fitting's radical of $\F_r$ by the free
$R$-module $T_s $ of the rank $s$, $\F_{r,s}=\F_r\oplus T_s$.}

\begin{lem} \label{cor:II-F-eqnoeth} \label{lem:221}
The free metabelian Lie algebra $\F_r$ is equationally Noetherian.
\end{lem}
\proof Let $F_r \left[X\right]$ be the free $F_r $-algebra
generated by the alphabet $X = \left\{ x_1 ,\ldots,x_n \right\} $.
Let $V_\M $ be the  verbal ideal of $F_r \left[X\right]$ that
defines the variety of all metabelian Lie algebras, $V_\M = \id
\left< (ab)(cd)\mid| {a,b,c,d \in F_r \left[X\right]} \right> $.
Consider the factor-algebra of the free $F_r$--algebra by $V_\M$:
$$
{\raise0.7ex\hbox{${F_r
[X]}$} \!\mathord{\left/
 {\vphantom {{F_r [X]} {V_\M }}}\right.\kern-\nulldelimiterspace}
\!\lower0.7ex\hbox{${V_\M}$}} \cong F_r  \ast_\M F_\M (X),
$$
where $\ast_\M$ stands for the free metbelian product of
metabelian Lie algebras and $F_\M (X) = F_n$ is the free
metabelian Lie algebra generated by the set $X$. Consequently, the
obtained factor-algebra is isomorphic to the free algebra $F_{r +
n}$. The Lie algebra  $F_{r + n}$ is metabelian, thus is Notherian
(in the usual classical sense, i.e. every its ideal is finitely
generated).

Fix an arbitrary system $S \subseteq F_r \left[X\right]$. To prove
the lemma it suffices to find a finite subsystem $S_0 \subseteq S$
so that $V_{F_r} (S) = V_{F_r} (S_0)$. Let $\bar S$ be the image
of  $S$ in the factor-algebra $F_r \ast_\M F_\M (X)$ and let  $I$
be the ideal generated by $\bar S$. Since $F_r  \ast_\M F_\M (X)
\cong F_{r + n} $ is a metabelian algebra, the ideal $I$ is
finitely generated. Choose a finite subsystem $\bar S_0 \subseteq
\bar S$ so that the set $\bar S_0 $ generates the ideal $I$.
Consider the pre-images of $S$ in $F_r [X]$ and take $S_0
\subseteq S$ to be an injective subset of pre-images for the set
$\bar S_0$. In the above notation, thanks to the choice of $S_0$,
it is clear that $V_{F_r } (S) = V_{F_r } (S_0 )$.\hfill
$\blacksquare$

\begin{rem}
The argument of Lemma \ref{lem:221} holds for an arbitrary
finitely generated metabelian Lie algebra $A$, i.e. every finitely
generated metabelian Lie algebra is equationally Noetherian.
\end{rem}

Let $\V$ be an arbitrary variety of Lie algebras and let $A \in
\V$. We shall make use of the following notation. Let $F_\V (X)$
be the free Lie algebra with the free base $X$ in the variety $\V$
and denote by $A_\V \left[ X \right] = A \ast_\V F_\V (X)$ the
$\V$-free product of $A$ and $F_\V (X)$ (for definitions see
\cite{Bah}). Suppose that $\varphi $ is the canonical homomorphism
from the free $A$--algebra $\AX$ to $A_\V [X]$. Then for every
system of equations $S \subseteq A\left[X\right]$ holds $\ker
\varphi \subseteq \rad _A (S)$. Which implies that
$$
\Gamma (S) = {\raise0.7ex\hbox{${\AX}$} \!\mathord{\left/
 {\vphantom {{\AX} {\rad_A (S)}}}\right.\kern-\nulldelimiterspace}
\!\lower0.7ex\hbox{${\rad_A (S)}$}} \cong _A
{\raise0.7ex\hbox{${A_\V \left[X\right]}$} \!\mathord{\left/
 {\vphantom {{A_\V \left[X\right]} {\overline {\rad_A (S)} }}}\right.\kern-\nulldelimiterspace}
\!\lower0.7ex\hbox{${\overline {\rad_A (S)} }$}},\quad \overline
{\rad_A (S)}  = \rad_A (\overline{S}).
$$
Where here ${\overline {\rad_A (S)} }$ denotes the image of
$\rad_A (S)$ in $A_\V [X]$ under $\varphi $. The definition of
$\rad_A (\overline{S})$ coincides with the one of $\rad_A (S)$
with  all the preliminary notions given with respect to the
algebra $A_\V [X]$

\begin{rem} \label{rem:II-tovar}
Let $\M$ be the variety of all metabelian Lie algebras and let
$A\in \M$. From the above discussion follows that all arguments
for the radicals of systems of equations over $A$ and respective
coordinate algebras can be performed in the metabelian Lie algebra
$A_{\M} \left[X\right]=A \ast_{\M} \F\left[X\right]$.
\end{rem}

Let $S$ be a system of equations over $\F_r $. On behalf of Remark
\ref{rem:II-tovar}, we may assume that $S \subseteq (\F_r)_{\M}
\left[X\right]=\F_r\ast_{\M} \F_r\left[X\right]$. Denote by $I_X =
\left< X \right> $ the ideal of $(\F_r )_{\M} \left[X\right]$,
generated by the alphabet $X$. Consider an equation $f(x_1
,\ldots,x_n ) \in S$. We next write $f$ as a sum of homogeneous
(by the variables $x_1, \ldots, x_n$) monomials:
$$
f = c + x_1 h_1  + \cdots + x_n h_n  + g(x_1 ,\ldots,x_n ),
$$
where here $c \in \F_r $, $h_i  \in R$ are polynomials, $i =
1,\ldots,n$, $g(x_1 ,\ldots,x_n ) \in I_X^2$. Note that the
expression $x_i h_i$ is not uniquely defined but the choice of
presentations of $h_1, \ldots h_n$ is not significant, see
\cite{AGFMLA1}. Since any solution of the system $S$ is a sum of
the linear part and the part from $\Fit (F_r)$, to solve the
system $S$ it is convenient to write the variables $x_i $, $i =
1,\ldots,n$ as a sum of two variables:
$$
  x_i  = z_i  + y_i ,\quad y_i  \in \Fit(\F_r ),\;z_i  = \alpha _{i1} a_1  + \cdots + \alpha _{ir} a_r
,\quad \alpha _{ij}  \in k.
$$
This increases the number of variables to $nr+n$: the variables
$\alpha _{ij}$, $i = 1,\ldots,n$, $j = 1,\ldots,r$ are valued in
$k$ and the variables $y_1 ,\ldots,y_n $ in $\Fit(\F_r)$.

However, this increase of the number of variables simplifies the
system $S$. Indeed, substituting the variables $x_i$ in the form
$z_i  + y_i$ into the equation $f = 0$ parts the system into two.
We separately equate to zero the linear part and the part from
$\Fit(\F_r)$. The first system is a regular linear system of
equations, which is solved using the methods of linear algebra. If
the correspondent system of equations over $k$ is inconsistent,
then so is the initial system $S$ over $\F_r$. Suppose that the
linear system of equations is consistent. To every its solution
$\alpha '_{ij}$, $i = 1,\ldots,n$, $j = 1,\ldots,r$ we associate
secondary module system of equations over $\Fit(\F_r)$ with
coefficients in $R$.

\section{Universal Axioms and the $\Phi_r$-Algebras}
\label{ss:II-univax-PhiR} \label{sec:3}

The above notation apply. In this section we formulate two
collections of seven series of universal  axioms $\Phi _r$ and
$\Phi '_r$ in the languages $L$ and $L_{\F_r}$. In the next
section we prove that $\Phi _r$ and $\Phi '_r$  axiomatise the
universal classes $\ucl(\F_r)$ and $\F_r-\ucl(\F_r)$. Most of
these formulas are the formulas of the first order language $L$.
Consequently they are the formulas  of both $\Phi _r $ and $\Phi
'_r $. We, therefore, write these series simultaneously, pointing
the differences
 between  $\Phi _r $ and $\Phi '_r $.

Since the algebra $\F_r $ is metabelian we write the metabelian
identity
$$
\Phi 1: \ \forall x_1 ,x_2 ,x_3 ,x_4 \quad (x_1 x_2 )(x_3 x_4 ) =
0.
$$

On account of Lemma 3.3 in \cite{AGFMLA1} the algebra $\F_r$
satisfies two following axioms
$$
\Phi 2: \ \forall x\forall y\quad xyx = 0 \wedge xyy = 0
\rightarrow xy = 0.
$$
$$
\Phi 3: \ \forall x\forall y\forall z\quad x \ne 0 \wedge xy = 0
\wedge xz = 0 \rightarrow yz = 0.
$$
The universal formula $\Phi 3$ is called the CT-axiom (commutative
transitivity axiom).

Let  $\N_2$ be the quasi variety of all Lie algebras defined by
$\Phi 1$ and $\Phi 2$  and let  $\N_3$ be the universal class
axiomatised by $\Phi 1$, $\Phi 2$ and $\Phi 3$.

\begin{lem} \label{lem:31}
Let $B \in \N_2 $. Then $B$ is a metabelian Lie algebra such that
{\rm $\Fit(B)$} and every nilpotent subalgebra of $B$ are Abelian.
\end{lem}
\proof By axiom $\Phi 1$ the algebra $B$ is metabelian. Let $C$ be
a nilpotent subalgebra of $B$ and $c_1 ,c_2  \in C$. Suppose that
$c_1  \circ c_2  \ne 0$. Then on account of Lemma 2.2 from
\cite{AGFMLA1} there exists a two-generated nilpotent subalgebra
$D =  \left< d_1 ,d_2  \right>$ of class two in $C$. It is
essentially immediate that $d_1 d_2 d_1  = d_1 d_2 d_2 = 0$, while
$d_1 d_2  \ne 0$, deriving a contradiction to $\Phi 2$. Finally,
recall that if every nilpotent subalgebra of $B$ is Abelian then
$\Fit(B)$ is Abelian (see Lemma 2.4 in \cite{AGFMLA1}). \hfill
$\blacksquare$

\begin{lem}\label{lem:32}
Let $B \in \N_3$ and let {\rm $a \in B \smallsetminus \Fit(B)$}.
Then for an arbitrary non-zero element $b$ from  {\rm $\Fit(B)$}
holds $ab \ne 0$.
\end{lem}
\proof Assume the converse: $ab = 0$ for a non-zero element $b \in
\Fit(B)$. For $\Fit(B)$ is Abelian, $bd = 0$ for any $d \in
\Fit(B)$. By the CT-axiom it follows that $ad = 0$, i. e. $a$
commutes with elements from  $\Fit(B)$. Consequently, see Lemma
2.4 in \cite{AGFMLA1}, $a \in \Fit(B)$. This contradicts the
assumption of the lemma. \hfill $\blacksquare$

We next introduce universal formula $\sFit(x)$ of the language $L$
with one variable $x$. The formula $\sFit(x)$ defines the
Fitting's radical
\begin{equation} \label{eq:fit1}
\sFit(x) \equiv (\forall y\;\,xyx = 0).
\end{equation}
The analogue of Formula (\ref{eq:fit1}) in the language $L_{\F_r
}$ is
\begin{equation} \label{eq:fit2}
\sFit '(x) \equiv (\,\mathop  \wedge \limits_i \,(xa_i x = 0)\,).
\end{equation}

\begin{lem} \label{lem:II-bN3} \label{lem:33}
Let $B \in \N_3$. Then the truth domain of Formula \/ {\rm
(\ref{eq:fit1})} is {\rm $\Fit(B)$}. In the event that $B$ is an
$\F_r$--algebra the truth domain of  {\rm $\sFit(x)$} is also {\rm
$\Fit(B)$}.
\end{lem}
\proof According to Lemma \ref{lem:31}, the Fitting's radical
$\Fit(B)$ is Abelian. It, therefore, is contained in the truth
domain of the Formula (\ref{eq:fit1}) (Formula (\ref{eq:fit2}),
respectively). Conversely, if $b \in B$ satisfies Formula
(\ref{eq:fit1}) then the ideal $I =  \left< b \right>$ is Abelian
and consequently $b \in \Fit(B)$.

For the case of $F_r$--algebras, we note that if $b \notin
\Fit(B)$ then, by axiom $\Phi 3$, for an element $a_i $ from the
free base of $F_r $ holds $ba_i  \ne 0$. Lemma \ref{lem:32},
therefore implies that $ba_i b \ne 0$. \hfill $\blacksquare$
\medskip
\newline
\textit{\textbf{Nota Bene}} \textit{We next restrict ourselves to
the case of a finite field $k$. In which case the vector space
 {\rm ${\raise0.7ex\hbox{${\F_r }$} \!\mathord{\left/
 {\vphantom {{\F_r } {\Fit(\F_r )}}}\right.\kern-\nulldelimiterspace}
\!\lower0.7ex\hbox{${\Fit(\F_r )}$}}$} is finite and its dimension
over $k$ is $r$.}
\medskip

\begin{lem} \label{lem:II-linindmodfit} \label{lem:34}
Let $k$ be a finite field and $n \in \mathbb{N}$, $n \le r$. Then
existential formula {\rm
$$
\varphi (x_1 ,\ldots,x_n) \equiv (\mathop \bigwedge
\limits_{(\alpha _1 ,\ldots,\alpha _n ) \ne \bar 0} \neg \sFit
(\alpha _1 x_1  + \cdots + \alpha _n x_n )).
$$}
of the language $L$ is true on the elements $\left\{ b_1
,\ldots,b_n \right\}$ of $\F_r$ if and only if \/ $b_1
,\ldots,b_n$ are linearly independent modulo {\rm $\Fit(\F_r )$}.
\end{lem}
\proof First we note that, since $k$ is finite, there exist only a
finite number of $n$-tuples  $\alpha _1 ,\ldots,\alpha _n \in k$.
In what it follows that $ \varphi (x_1 ,\ldots,x_n ) $ is a
formula of the language $L$. Now, since $\Fit(x)$ is a universal
formula, the negation  $\neg \Fit(x)$ is an existential formula
and so is  the formula $\varphi (x_1 ,\ldots,x_n )$.

Let $\left\{ b_1 ,\ldots,b_n \right\}$ be a system of elements
from the truth domain of the formula $\varphi $. By Lemma
\ref{lem:33} this system is linearly independent modulo
$\Fit(F_r)$. The converse is obviously also true.

Resulting from Lemma \ref{lem:II-bN3}, the truth domain of
$\varphi (x_1 ,\ldots,x_n)$ in the algebra $B \in \N_3$ is  the
set of all linearly independent modulo $\Fit(B)$ $n$-tuples, $n
\le r$. Therefore, the formula $\varphi$ formalises the notion of
linear independence modulo the Fitting's radical for a tuple of
elements, provided that the field $k$ is finite.  The formula
$\varphi$ is very convenient to use in the language $L$. In the
language $L_{\F_r}$ there is another, more simple way to test
whether the elements $b_1 ,\ldots,b_n $ are linearly independent
modulo $\Fit(B)$.

\begin{lem} \label{lem:35}
Let $B$ be an $\F_r$--Lie algebra, $B \in \N_3$ and let $c_1
,\ldots,c_n$, $n \le r$ be linearly independent elements modulo
the Fitting's radical of the designated copy of $\F_r$. Then the
elements $c_1 ,\ldots,c_n$ are linearly independent modulo {\rm
$\Fit(B)$}.
\end{lem}
\proof Suppose that a non-trivial linear combination $c = \alpha
_1 c_1 + \cdots + \alpha _n c_n$, $\alpha _i  \in k$ lies in
$\Fit(B)$. Since $\Fit (B)$ is Abelian, $a_1 a_2 c = 0$. This
derives a contradiction, for $F_r $ is a $U$-algebra. \hfill
$\blacksquare$

With the help of the formula $ \varphi$ we next write the
dimension axiom
$$
\Phi 4: \forall x_1 ,\ldots,x_{r + 1} \;\neg \varphi (x_1
,\ldots,x_{r + 1} ).
$$
For $\varphi $ is an existential formula, the formula $\neg
\varphi $ is a universal formula, thus so is the formula $\Phi 4$.
This axiom postulates that the dimension of the factor-space
${\raise0.7ex\hbox{${B}$} \!\mathord{\left/
 {\vphantom {{B } {\Fit(B )}}}\right.\kern-\nulldelimiterspace}
\!\lower0.7ex\hbox{${\Fit(B)}$}}$ is lower than or equals $r$,
provided that $B \in \N_3 $.

Recall that  $\Fit(\F_r )$ allows the structure of a module over
the ring of polynomials $R = k\left[x_1 ,\ldots,x_r \right]$. The
series of axioms $\Phi 5$, $\Phi' 5$,  $\Phi 6$, $\Phi 7$ and
$\Phi' 7$ express module properties of  $\Fit(\F_r )$. We use
module notation here, i.e. the multiplication of elements of an
algebra on elements from $R$. By this notation (we refer to
\cite{AGFMLA1} for details) we mean that the polynomial $f(x_1
,\ldots,x_n)$, $n \le r$ rewrites into the signature of metabelian
Lie algebras.

The Fitting's radical of the free metabelian Lie algebra is a
torsion free module  over the ring $R$. Consequently we write the
following infinite series of axioms. For every non-zero polynomial
$f \in k \left[x_1 ,\ldots,x_n \right]$, $n \le r$ write
\begin{gather} \notag
\begin{split}
\Phi 5: \ \forall z_1 ,z_2 \;\forall x_1 ,\ldots,x_n \quad (z_1
z_2 \cdot f(x_1 ,\ldots,x_n ) & = 0 \, \wedge \, z_1 z_2  \ne 0)
\rightarrow \\
& \rightarrow (\neg \, \varphi (x_1 ,\ldots,x_n )).
\end{split}
\end{gather}
Since  $\varphi (x_1 ,\ldots ,x_n )$ is a $\exists$-formula the
formula $\Phi 5$ is a $\forall$-formula.

For the collection $\Phi '_r$, i.e. for axioms in the language
$L_{\F_r}$, this fact can be expressed in a more simple way,
$$
\Phi' 5: \ \forall z_1 ,z_2 \quad (z_1 z_2 \cdot f(a_1 ,\ldots,a_r
) = 0 \rightarrow z_1 z_2  = 0).
$$
Here $f$ is a non-zero polynomial from $k \left[x_1 ,\ldots,x_r
\right]$.

The main advantage of this formula is that it does not involve the
formula $\varphi$, which implies that the restriction on the
cardinality of the field $k$ is not significant.

Let $\N_5 $ and $\N'_5 $ correspondingly be the universal classes
generated by the series of axioms  $\Phi 1-\Phi 5$ and $\Phi
1-\Phi 4, \Phi' 5$.

\begin{lem} \label{lem:36} \
The class $\N_5$ is the class of all $U$-algebras $B$ such that \/
{\rm $\dim {\raise0.7ex\hbox{${B}$} \!\mathord{\left/
 {\vphantom {{B } {\Fit(B )}}}\right.\kern-\nulldelimiterspace}
\!\lower0.7ex\hbox{${\Fit(B)}$}} \le r$}.
\newline
The class  $\N'_5 $ is the class of all $\F_r$--\/$U$-algebras $B$
such that \/ {\rm $\dim {\raise0.7ex\hbox{${B}$} \!\mathord{\left/
 {\vphantom {{B } {\Fit(B )}}}\right.\kern-\nulldelimiterspace}
\!\lower0.7ex\hbox{${\Fit(B)}$}} = r$}.
\end{lem}
\proof Let $B$ be a non-abelian Lie algebra from the class $\N_5$.
Then, according to Lemma \ref{lem:31}, $\Fit(B)$ is an Abelian
ideal. Assume that $\dim {\raise0.7ex\hbox{$B$} \!\mathord{\left/
{\vphantom {B {\Fit(B)}}}\right.\kern-\nulldelimiterspace}
\!\lower0.7ex\hbox{${\Fit(B)}$}} = n$. The inequality $n \le r$
follows immediately from Lemma \ref{lem:34} and axiom $\Phi 4$. In
the event that  $B$ is an $F_r$--algebra  $n$ equals $r$ (see
Lemma \ref{lem:35}).

We next show that $B$ is a  $U$-algebra. The Fitting's radiacal
$\Fit(B)$ admits a structure of a module over the ring $k\left[x_1
,\ldots,x_n \right]$. We show that $\Fit(B)$ is torsion-free. Take
an element $0 \ne b \in \Fit(B)$ and a non-zero polynomial $f(x_1
,\ldots,x_n )$. By Lemma \ref{lem:32}, $ba \ne 0$ for any  $a \in
B\smallsetminus \Fit(B)$. By axiom $\Phi 4$,  $(ba) \cdot f(x_1
,\ldots,x_n ) \ne 0$. Since $b \in \Fit(B)$, the product $ba$ can
be written as: $ba = b \cdot g(x_1 ,\ldots,x_n )$, where $g(x_1
,\ldots,x_n )$ is a linear non-zero polynomial. In what follows
that $b \cdot f(x_1 ,\ldots,x_n ) \ne 0$, that $\Fit(B)$ is a
torsion-free module, and therefore $B$ is a $U$-algebra.

Conversely, let $B$ be a metabelian $U$-algebra and let $\dim
{\raise0.7ex\hbox{$B$} \!\mathord{\left/ {\vphantom {B
{\Fit(B)}}}\right.\kern-\nulldelimiterspace}
\!\lower0.7ex\hbox{${\Fit(B)}$}} = n$ , $n \le r$. We show that
$B$ lies in  $\N_5 $. From elementary properties of metabelian
$U$-algebras (see Theorem 3.4 \cite{AGFMLA1}) it follows that $B$
lies in  $\N_3 $. From Lemmas \ref{lem:34} and \ref{lem:35} it is
immediate that $B$ lies in  $\N_5$ ($\N'_5 $). \hfill
$\blacksquare$

According to  Corollary 2.4 in \cite{AGFMLA1}, every $n$-tuple of
linearly independent modulo $\Fit (\F_r)$ elements $\left\{ b_1
,\ldots,b_n \right\}$, $n \le r$  of  $\F_r$ freely generates the
free metabelian Lie algebra of the rank $n$. In the event that $k$
is finite, the statement of Corollary 2.4 for algebras from $\N_5$
and $\N'_5 $ can be written by the means of universal formulas.

For every non-zero Lie polynomial $l(a_1 ,\ldots,a_n )$, $n \le r$
of the letters $a_1 ,\ldots,a_r $ from the free base of $\F_r$
write
$$
\Phi 6: \ \forall x_1 ,\ldots,x_n \quad \varphi (x_1 ,\ldots,x_n )
\rightarrow (l(x_1 ,\ldots,x_n ) \ne 0).
$$
Since  $\varphi (x_1 ,\ldots ,x_n )$ is a $\exists$-formula the
formula $\Phi 6$ is a $\forall$-formula.

Denote by $\N_6 $
 and $\N'_6 $, correspondingly, the universal classes
generated by the series of axioms $\Phi 1-\Phi 6$ and $\Phi 1-\Phi
4, \Phi' 5, \Phi 6$.

\begin{lem} \label{lem:II-N6} \label{lem:37}
Let $B = \F_n  \oplus M$ be the direct module extension of $\F_n$
(see Section 4.3 in \cite{AGFMLA1}), where here $M$ is a torsion
free module over $k\left[x_1 ,\ldots,x_n \right]$, $n \le r$. Then
$B\in \N_6$ and, in the event that $n = r$, $B \in \N'_6 $.
\end{lem}
\proof Note that, according to Lemma \ref{lem:36}, the free
metabelian Lie algebra $F_n $ of the rank $n \le r$ lies in
$\N_5$. Therefore, it is clear that $F_n \in \N_6$. Since
$\ucl(F_n ) = \ucl(B)$ (see Proposition 4.4 in \cite{AGFMLA1}),
the algebra $B$ also lies in $\N_6 $. \hfill $\blacksquare$

Unfortunately not every finitely generated algebra from $\N_6$ has
the form $\F_n  \oplus M$. On the other hand, all finitely
generated algebras from  $\ucl (\F_r)$ and $\F_r -\ucl (\F_r)$
have this form (see Theorems \ref{thm:II-ucl} and
\ref{thm:II-F-ucl}). The point is that every Lie algebra from
$\N_6$ (correspondingly, $\N'_6 $) is obtained from  $\F_n$, $n
\le r$ (correspondingly from $\F_r $) by the means of an extension
of its Fitting's radical.  But in general, this extension is not
the direct module extension, i.e. the Fitting's radical of the
algebra from  $\N_6$ (or from  $\N'_6$) is not the \emph{direct}
sum of a new module and the initial Fitting's radical. To narrow
the classes $\N_6 $ and $\N'_6 $ we need to write the final the
most sophisticated series of axioms $\Phi 7$ and $\Phi' 7$.

We first introduce  higher-dimensional analogues of Formulas
(\ref{eq:fit1}) and (\ref{eq:fit2})
\begin{equation} \label{eq:fitm1}
\sFit(y_1,\ldots,y_l ;\,x_1 ,\ldots,x_n ) \equiv (\,\mathop
\bigwedge \limits_{i = 1}^n \,(y_1 x_i y_1 = 0)\,) \wedge \ldots
\wedge (\,\mathop  \bigwedge \limits_{i = 1}^n \,(y_l x_i y_l  =
0)),
\end{equation}
\begin{equation} \label{eq:fitm2}
    \sFit'(y_1 ,\ldots,y_l ) \equiv \sFit'(y_1 ) \wedge \ldots \wedge \sFit'(y_l ).
\end{equation}

 Formula (\ref{eq:fitm2}) defines the subset
$$
\underbrace {\Fit(B) \times \cdots \times \Fit(B)}_l \subseteq
\underbrace {B \times \cdots \times B}_l,
$$
provided that $B$ is an $\F_r$--Lie algebra from $\N_3$.

Let $B$ be an algebra from $\N_6$ and let $\left\{ b_1 ,\ldots,b_n
\right\}$ be a linearly independent modulo $\Fit(B)$ set of
elements from $B$, $n \le r$, $n \ge 2$.
    The truth domain of the formula $\sFit(y_1 ,\ldots,y_l ;\,b_1 ,\ldots,b_n)$ is
    $\underbrace {\Fit(B) \times \cdots \times \Fit(B)}_l$.
In the event that $n = 1$ the same holds, provided that
    $\dim {\raise0.7ex\hbox{$B$}
\!\mathord{\left/
 {\vphantom {B {\Fit(B)}}}\right.\kern-\nulldelimiterspace}
\!\lower0.7ex\hbox{${\Fit(B)}$}} = 1$.

Before we introduce the final series of axioms we need to explain
the syntaxis of these formulas. We begin with the series of axioms
$\Phi' 7$ in the language $L_{\F_r }$. Let $S$ be a fixed finite
system of module equations with variables $y_1 ,\ldots,y_l $ over
the module $\Fit(\F_r )$. Every equation from $S$ has the form
$$
h = y_1 f_1 (\bar x) + \cdots + y_l f_l (\bar x) - c = 0,\;\quad c
= c(a_1 ,\ldots,a_r ) \in \Fit(\F_r ),
$$
where here $\bar x = \{ x_1 ,\ldots,x_r \} $ is a vector of
variables and $f_1 ,\ldots,f_l  \in R = k\left[\bar x\right]$.
Suppose that $S$ is inconsistent over  $\Fit(\F_r )$. This fact
can be  easily written in the signature of a module. The system
$S$ gives rise to a system of equations $S_1$ over $\F_r$. Replace
every
 module equation $h_i = 0$ from $S$ by the equation $h'_i  = 0$, $i =
1,\ldots,m$ in the signature of $L_{\F_r}$ (see \cite{AGFMLA1}).
This results a system of equations $S_1$ over $\F_r$. By every
inconsistent module system of equations $S$ write
$$
\Phi '7: \ \psi '_S  \equiv \forall y_1 ,\ldots,y_l \quad
\sFit'(y_1 ,\ldots,y_l ) \rightarrow \mathop  \bigvee \limits_{i =
1}^m \;h'_i (y_1 ,\ldots,y_l ) \ne 0.
$$
Notice that the restriction on the cardinality of the field $k$ is
not used.

Denote by $\N'_7 $ the universal class  axiomatised by  $\Phi
1-\Phi 4, \Phi' 5, \Phi 6, \Phi' 7$.

\begin{lem} \label{lem:II-directdecomp} \label{lem:38}
Let  $B$ be a finitely generated $\F_r $--algebra  from $\N'_7$.
Then $B$ is $\F_r $--isomorphic to the algebra $\F_r  \oplus M$
for some finitely generated torsion free module  $M$ over $R$.
\end{lem}
\proof Suppose that the elements $a_1 ,\ldots,a_r $ generate the
designated copy of $F_r$ in $B$. By Lemma \ref{lem:35} these
elements are linearly independent modulo $\Fit(B)$. For  $B$ is a
finitely generated algebra from  $\N'_7 $, from Lemma \ref{lem:36}
we conclude that $\Fit(B)$ is a finitely generated torsion-free
module over the ring $R$ (see Lemma 2.6 in \cite{AGFMLA1}). It,
therefore, suffices to prove that the submodule $\Fit(F_r )$ of
the module $\Fit(B)$ is a direct summand.

For the sake of brevity, set $P = \Fit(B)$, $N = \Fit(F_r )$. We
prove that the factor-module ${\raise0.7ex\hbox{$P$}
\!\mathord{\left/ {\vphantom
{PN}}\right.\kern-\nulldelimiterspace} \!\lower0.7ex\hbox{$N$}}$
is torsion-free. Assume the converse, then there exists an element
$d \in \Fit(F_r )$ and a non-zero polynomial $f(x_1 ,\ldots,x_r )$
so that the equation
\begin{equation} \label{eq:eqmod}
y \cdot f(x_1 ,\ldots,x_r ) = d
\end{equation}
is compatible over $\Fit(B)$ (let $y_1 \in \Fit(B)$ be its
solution) but incompatible over $\Fit(F_r )$. Assume that, on the
contrary, there exists a solution $y_2 \in \Fit(F_r )$. In what
follows that $(y_1-y_2) \cdot f(x_1 ,\ldots,x_r )=0$, deriving a
contradiction with the fact that $P$ is torsion-free.

Equation (\ref{eq:eqmod}) gives rise to an incompatible module
system of equations $S_0 $ over  $\Fit(F_r )$. In what it follows
that the axiom $\psi '_{S_0 }$ from the series  $\Phi' 7$ is true
in $B$, which is obviously false. Further, if
${\raise0.7ex\hbox{$P$} \!\mathord{\left/   {\vphantom
{PN}}\right.\kern-\nulldelimiterspace} \!\lower0.7ex\hbox{$N$}}$
is a free module then $P = N \oplus M$ and $B = F_r  \oplus M$. So
we assume that  ${\raise0.7ex\hbox{$P$} \!\mathord{\left/
{\vphantom {P N}}\right.\kern-\nulldelimiterspace}
\!\lower0.7ex\hbox{$N$}}$ is not a free module and that the images
of the elements $m_1 ,\ldots,m_l  \in P$ are its generators. Let
$S$ be a finite system of module relations for the generators
$\bar m_1 ,\ldots,\bar m_l$ of the module ${\raise0.7ex\hbox{$P$}
\!\mathord{\left/ {\vphantom {P
N}}\right.\kern-\nulldelimiterspace} \!\lower0.7ex\hbox{$N$}}$. If
in the left-hand sides of the relations we replace the generators
by the variables $y_i$'s and in the right-hand side we write the
values of the left-hand side relations with $y_i$'s valued as
$m_i$'s, we obtain a finite system $S_1$ of module equations over
$N$. The system $S_1 $ has a solution $\left\{ m_1 ,\ldots,m_l
\right\}$ in $P$. Consequently, by the axioms of the series $\Phi'
7$, the system $S_1$ has a solution $\left\{ c_1 ,\ldots,c_l
\right\} $ in  $\Fit(F_r )$. Denote by $M$ the submodule of the
module $P$ generated by the elements $m_i  - c_i  = m'_i$, $i =
1,\ldots,l$. Clearly, $M$ is isomorphic to the module
${\raise0.7ex\hbox{$P$} \!\mathord{\left/ {\vphantom {P
N}}\right.\kern-\nulldelimiterspace} \!\lower0.7ex\hbox{$N$}}$,
and consequently  $P = N \oplus M$. In what follows that $B = F_r
\oplus M$ (see Lemma 4.8 \cite{AGFMLA1}). \hfill $\blacksquare$

Now we turn to the collection $\Phi_r$ in the language $L$. Let
$B\in \N_6$ and let $C$ be its subalgebra generated by the
elements $c_1,\ldots,c_n$, $n \le r$. Assume that $c_1,\ldots,c_n$
are linearly independent modulo $\Fit(B)$. Then $C$ is isomorphic
to $\F_n $.  Although, even in the event that $B=\F_r$ the
subalgebra
 $C\varsubsetneq \F_r$ does not yield to the decomposition $\F_r  \oplus M$.
To avoid this problem we use the notion of $\Delta $-localisation
of $\F_r $ (see Section 4.2 in \cite{AGFMLA1}).

The basic idea of the formula $\Phi 7$ is to write it in such a
way that an analogue of Lemma \ref{lem:II-directdecomp} holds for
$\Delta $-local Lie algebras from $\N_7$.

We next introduce some auxiliary notation. Let  $S$ be a finite
system of module equations with variables $y_1 ,\ldots,y_l$ over
$\Fit(\F_n )$, $n \le r$  and let $f(x_1,\ldots,x_n )$ be a
polynomial  from $R \smallsetminus \Delta$, where here $\Delta=
\left< x_1 ,\ldots , x_r \right>$.  Denote by $UD(f)$ the
collection of all unitary divisors of $f$, $UD(f) \subseteq R
\smallsetminus \Delta $. For each $f$ and each
 $\alpha  \in \underbrace {UD(f) \times \cdots \times UD(f)}_l$,
 $\alpha  = (d_1 ,\ldots,d_l )$ define a system of equations
 $S_{f,\alpha }$. The system $S_{f,\alpha }$ is obtained from $S$
by multiplying the equations from $S$ by the polynomial $d = d_1
\cdots d_l $ and dividing the coefficient of the term $y_i $ by
$d_i $.

\begin{lem} \label{lem:39}
Let $B$
be an algebra from $\N_6$. In the above notation, the system $S$
is consistent over {\rm $\Fit_\Delta  (B)$} if and only if there
exist $f(x_1 ,\ldots,x_n ) \in R\smallsetminus \Delta $ and
$\alpha \in UD(f) \times \cdots \times UD(f)$ such that the system
$S_{f,\alpha}$ is consistent over {\rm $\Fit(B)$}.
\end{lem}
\proof The proof is straightforward. \hfill $\blacksquare$

\begin{rem}
If $S$ is inconsistent over {\rm $\Fit_\Delta  (B)$}, then for all
$f$ and $\alpha $ the system $S_{f,\alpha }$ is inconsistent over
{\rm$\Fit_\Delta  (B)$}.
\end{rem}

Suppose that  $S$ is inconsistent over the Fitting's radical
{\rm$\Fit_\Delta (\F_n )$} of $\Delta $-local Lie algebra $(\F_n
)_\Delta$ system of $m$ module equations over $\Fit(\F_n )$. For
every $n \in \mathbb{N}$, $n \le r$ and a system $S$ we write
\begin{gather} \notag
\begin{split}
\Phi 7: \ \psi _{n,S}  \equiv \forall x_1 ,\ldots,x_n \;\forall
y_1 ,\ldots,y_l & \quad \varphi (x_1 ,\ldots,x_n ) \wedge
\sFit(y_1 ,\ldots,y_l ) \rightarrow \\
& \rightarrow \mathop \bigvee \limits_{i = 1}^m h_i (y_1
,\ldots,y_l ;\,x_1 ,\ldots,x_n ) \ne 0.
\end{split}
\end{gather}
The Lie polynomials $h_i $, $i = 1,\ldots,m$ are obtained from the
system $S$. Consider the $i$-th polynomial from $S$. It has the
form
\begin{gather} \notag
\begin{split}
h'_i = y_1 f_1 (x_1 ,\ldots,x_n ) + \cdots + & y_l f_l (x_1
,\ldots,x_n) - c = 0,\\
& f_i  \in R,\;c = c(a_1 ,\ldots,a_n ) \in \Fit(\F_n ).
\end{split}
\end{gather}
Its interpretation in the signature of  Lie algebras (see
\cite{AGFMLA1}) with every occurrence of $a_j $ in $c(a_1
,\ldots,a_n )$ replaced by $x_j $, $j = 1,\ldots,n$ results the
polynomial $h_i$.

\begin{lem} \label{lem:310}
The axioms of the series $\Phi 7$ are true in the free metabelian
Lie algebra $\F_r$.
\end{lem}
\proof Let $S$ be a finite incompatible over $\Fit_\Delta  (F_n )$
system of module equations over $\Fit(F_n )$, $n \le r$. In which
case it is clear that  $S$ is incompatible over $\Fit_\Delta
(F_r)$. Take a tuple of arbitrary linearly independent modulo
$\Fit(F_r )$ elements $c_1 ,\ldots,c_r \in F_r $. Denote by $C$
the subalgebra of $F_r $ generated by these elements. The
subalgebra $C$ is isomorphic to the algebra $F_r $. Therefore, the
system $S$ is incompatible over  $\Fit_\Delta  (C)$ and thus $h_i
(b_1 ,\ldots,b_l ;\,c_1 ,\ldots,c_n ) \ne 0$, $i = 1,\ldots,m$ for
any $b_1 ,\ldots,b_l  \in \Fit(C_\Delta  )$. But if we treat $C$
as a set it might not coincide with $F_r $. However we have
$C_\Delta = (F_r )_\Delta  $ and  $\Fit(C_\Delta  ) = \Fit_\Delta
(F_r )$ (see Proposition 4.2 \cite{AGFMLA1}). We, therefore,
obtain that the correspondent formula $\psi _{n,S}$ is true in the
algebra $F_r $. \hfill $\blacksquare$

Denote by $\N_7 $ the universal class  axiomatised by  $\Phi 1-
\Phi 7$.

\begin{lem} \label{lem:II-311} \label{lem:311}
If an algebra $B$ lies in $\N_7$ then its $\Delta $-localisation
$B_\Delta $ lies in $\N_7 $ and has the form $B_\Delta = (\F_n
)_\Delta   \oplus M_\Delta  $ for some $n \le r$ and some finitely
generated torsion free module $M$ over the ring $k\left[x_1
,\ldots,x_n \right]$. Furthermore, there exists an integer $s \in
\mathbb{N}$ such that the algebra $B$ is a subalgebra of
$\F_{r,s}$.
\end{lem}
\proof Let  $\dim {\raise0.7ex\hbox{$B$} \!\mathord{\left/
{\vphantom {B {\Fit(B)}}}\right.\kern-\nulldelimiterspace}
\!\lower0.7ex\hbox{${\Fit(B)}$}} = n$. By Lemma \ref{lem:36} $n
\le r$ and $B$ is a $U$-algebra. By the axioms of the series $\Phi
6$, a tuple of linearly independent modulo $\Fit(B)$ elements $b_1
,\ldots,b_n$ generates a subalgebra of $B$ isomorphic to $F_n$.

The fact that $B_\Delta   \in \N_7 $ is implied by the coincidence
of universal closures $\ucl(B) = \ucl(B_\Delta  )$, which is true
for an arbitrary  $U$-algebra $B$ (see Proposition 4.1
\cite{AGFMLA1}).

We  next show that if a finite system of module equations $S$ has
a solution in $\Fit_\Delta  (B)$, then it has a solution in
$\Fit_\Delta  (F_n )$. If $S$ is compatible over $\Fit_\Delta (B)$
then by Lemma \ref{lem:39} for some $f(x_1 ,\ldots,x_n ) \in
R\smallsetminus \Delta $ and some $\alpha  \in D(f) \times \cdots
\times D(f)$ the system $S_{f,\alpha } $ is compatible over
$\Fit(B)$. From the series of axioms $\Phi 7$ it, therefore,
follows that $S_{f,\alpha } $ has a solution in $\Fit_\Delta (F_n
)$ and consequently, on account of Lemma \ref{lem:39}, the system
$S$ has a solution in  $\Fit_\Delta  (F_n)$. To prove that
$B_\Delta$ has the form $(F_n )_\Delta   \oplus M_\Delta$ it
suffices to conduct an argument analogous to the one of Lemma
\ref{lem:38}. We leave this to the reader.

Since $(F_n )_\Delta \oplus M_\Delta = (F_n  \oplus M)_\Delta$
(see Lemma 4.5 in \cite{AGFMLA1}) and since $B$ embeds into
$B_\Delta$, the algebra $B$ is a finitely generated subalgebra of
a $\Delta$-local algebra $(F_n  \oplus M)_\Delta$. In what follows
that  $B$ embeds into the algebra $F_n  \oplus M$ (see Lemma 4.1
in \cite{AGFMLA1}). The module  $M$ embeds into the free module
$T_s$ of the rank  $s$ over $R$ (see \cite{Bourb}, \cite{Leng}).
This embedding gives rise to an embedding of the algebra $F_n
\oplus M$ into the algebra $F_{n,s} $ (see Lemma 4.8
\cite{AGFMLA1}), which is a subalgebra of $F_{r,s} $. We,
therefore, have proven that  $B$ embeds into the algebra
$F_{r,s}$. \hfill $\blacksquare$

Denote, correspondingly, by $\Phi _r $ and by $\Phi '_r $ the
universal classes axiomatised by $\Phi 1- \Phi 7$ and by $\Phi 1-
\Phi 4, \Phi' 5, \Phi 6, \Phi'7$ respectively.

\begin{defn}
The algebras from $\Phi _r $ and $\Phi '_r $ are termed,
correspondingly, $\Phi _r $-algebras and $\Phi '_r $-algebras.
\end{defn}

\begin{lem} \label{lem:II-3.12} \label{lem:312}
For every $r \in \mathbb{N}$ and every $n \le r$ and $m > r$ the
algebra $\F_n $ is a $\Phi _r $-algebra, while $\F_m $ is not a
 $\Phi _r $-algebra.
The algebra $\F_r$ lies in $\Phi '_r$ and  $\F_r \notin \Phi '_n
$, provided that $n \ne r$.
\end{lem}
\proof From the axioms $\Phi _r$ and the properties of $F_r$
proved above, the algebra $F_n$ is a $\Phi _r$-algebra. However,
the algebra, $F_m $, $m > r$ does not lie in $\Phi _r $, for the
dimension axiom is false if $m > r$. \hfill $\blacksquare$

\begin{cor} The following sequence of strict inclusions holds
$$
\Phi _1  \varsubsetneq \Phi _2 \varsubsetneq \ldots  \varsubsetneq
\Phi _r  \varsubsetneq \ldots
$$
The classes $\Phi '_{r_1}$ and $\Phi '_{r_2}$ are disjoint.
\end{cor}

\begin{cor} Let  $A$ be a
$\Phi _r $-algebra and let {\rm $\dim {\raise0.7ex\hbox{$A$}
\!\mathord{\left/
 {\vphantom {A {\Fit(A)}}}\right.\kern-\nulldelimiterspace}
\!\lower0.7ex\hbox{${\Fit(A)}$}} = n < r$}. Then $A$ is a $\Phi
_n$-algebra but not a $\Phi _m $-algebra, where here $m < n$.
\end{cor}

Finally, we emphasise the aspects which impose the restriction on
the cardinality of the field $k$.

First of all, the formula $\varphi$ is finite only because the
ground field $k$ is finite. Recall that the formula $\varphi$ is
written in the language $L$ and formalises the notion of linear
independence modulo the Fitting's radical. The formula $\varphi $
is used in all the axioms of the series $\Phi 4 - \Phi 7$.

In the enriched language $L_{\F_r}$ the formula $\varphi $ is
almost unnecessary. It  is  involved in neither $\Phi '5$ nor
$\Phi' 7$ and is  only used in the dimension axiom $\Phi 4$ and in
the series $ \Phi 6$. But the axioms of $ \Phi 6$ can be excluded
from the collection $\Phi '_r$.  The series $\Phi 6$ postulates
that every set of linearly independent modulo the Fitting's
radical elements $\left\{ b_1 ,\ldots,b_n \right\} $, $n \le r$
freely generates an algebra isomorphic to $\F_n $. This property
is important only in the proofs of Lemmas
\ref{lem:II-directdecomp} and \ref{lem:II-311}, which require only
the existence of such $n$-tuples. Since every $\F_r$--algebra
contains a designated copy of $\F_r$, which possesses such an
$n$-tuple, the series $\Phi 6$ can be omitted in the case of
$L_{\F_r}$.

However, the axiom $\Phi 4$ is significant. We can not write this
axiom without the use of the formula $\varphi $ and can not
exclude $\Phi 4$ from $\Phi '_r $. In the event that the field $k$
is finite the dimension axiom is true in  $\F_r $ and is true in
every algebra from $\F_r - \ucl(\F_r)$. But in the event that the
main field is infinite the algebra $\F_r $, $r \ge 2$ is
discriminated by $\F_2 $. Which essentially implies that $\F_r -
\ucl(\F_r)$ contains algebras $A$ of `unlimited' dimension of
${\raise0.7ex\hbox{$A$} \!\mathord{\left/
 {\vphantom {A {\Fit(A)}}}\right.\kern-\nulldelimiterspace}
\!\lower0.7ex\hbox{${\Fit(A).}$}}$

\section{Theorems on Universal Closures of the Algebra  $F_r$}
\label{sec:4}

In this section we formulate and prove several theorems on
universal classes $U_r  = \ucl(F_r )$ and $U'_r  = F_r - \ucl(F_r
)$. The main results here are
\begin{itemize}
    \item The set of axioms $\Phi_r$ ($\Phi'_r$)
          axiomatise the universal closure of the free
          metabelian Lie algebra of finite rank $r \ge 2$ over a finite
          field $k$,
    \item Given a description of the structure of finitely
          generated algebras from  $\F_r- \ucl (\F_r) $ and $\ucl
          (\F_r)$,
    \item Investigated the structure of irreducible algebraic sets over $\F_r $
          and the structure of correspondent coordinate algebras,
    \item Proved that the universal theory of the free metabelian Lie algebra is decidable
    in both $L$ and $L_{\F}$.
\end{itemize}

\subsection{Formulation of Main Results} \label{sec:41}

\begin{thm} \label{thm:II-ucl} \label{thm:41}
Let $A$ be an arbitrary finitely generated metabelian Lie algebra
over a finite field $k$. Then the following conditions are
equivalent
\begin{itemize}
    \item  {\rm $A \in \ucl(\F_r )$};
    \item  $A$ is a $\Phi _r$-algebra;
    \item  there exists $s \in \mathbb{N}$ such that
    $A$ is a subalgebra of $\F_{r,s}$.
\end{itemize}
\end{thm}

\begin{cor} \label{cor:II-ucl} \label{cor:41}
The universal closure {\rm  $\ucl(\F_r )$} of the free metabelian
Lie algebra  $\F_r$ is axiomatised by $\Phi _r $.
\end{cor}

\begin{thm} \label{thm:II-F-ucl} \label{thm:42}
Let $A$ be an arbitrary finitely generated metabelian $\F_r$--Lie
algebra over a finite field $k$. Then the following conditions are
equivalent
\begin{itemize}
    \item {\rm $A \in \F_r-\ucl(\F_r )$};
    \item  $A$ is a $\Phi' _r$-algebra;
    \item  $A$ is $\F_r$--isomorphic to the algebra
     $\F_{r}\oplus M$, where here $M$ is a torsion free module over
     the ring of polynomials $R=k\left[x_1 ,\ldots,x_r \right]$;
\end{itemize}
\end{thm}

\begin{cor} \label{cor:II-F-ucl} \label{cor:42}
The universal closure  {\rm $\F_r-\ucl(\F_r )$} of the free
metabelian Lie algebra $\F_r$ is axiomatised by $\Phi' _r $.
\end{cor}

Recall that the fact that the field $k$ is finite, is significant
for the axiom $\Phi 4$, for series of axioms $\Phi 5$, $\Phi 6$,
$\Phi 7$.

Two next theorems treat the decidability of universal theory of
the algebra  $F_r $ in the languages  $L$ and $L_{F_r }$.

\begin{thm}\label{thm:II-Phi} \label{thm:43}
Axioms $\Phi _r $ form a recursive set and the universal theory in
the language $L$ of the algebra $\F_r$ over a finite field $k$ is
decidable.
\end{thm}

\begin{thm} \label{thm:II-Phi'} \label{thm:44}
Axioms $\Phi' _r $ form a recursive set and the universal theory
in the language $L_{\F_r}$ of the algebra $\F_r$ (treated as a
$\F_r$--algebra) over a finite field $k$ is decidable.
\end{thm}

\begin{thm} \label{thm:45}
Compatibility problem for system of equations over the free
metabelian Lie algebra $\F_r$ is decidable.
\end{thm}

This result contrasts  with a result of V.A. Roman'kov on the
compatibility problem over some metabelian algebraic systems. In
\cite{Roman} he proves that this problem is undecidable for free
metabelian groups of a large enough rank. The argument of
\cite{Roman} holds for free metabelian Lie rings and for free
metabelian Lie algebras, provided that the compatibility problem
for the ground field is undecidable.

\subsection{Proofs of the Theorems} \label{sec:42}

\emph{\textbf{\underline{Proof of Theorem \ref{thm:41}}}}

\underline{$1 \rightarrow 2$} On behalf of Lemma
\ref{lem:II-3.12}, $\F_r$ is a $\Phi _r$-algebra. Consequently,
every (unnecessarily finitely generated) Lie algebra from
$\ucl(\F_r )$ is a $\Phi _r$-algebra.

\underline{$2 \rightarrow 3$} On account of Lemma
\ref{lem:II-311}, every finitely generated Lie algebra  from $\Phi
_r$ is a subalgebra of $\F_{r,s}$ for some $s \in \mathbb{N}$.

\underline{$3 \rightarrow 1$} Finally, from Proposition 4.4 in
\cite{AGFMLA1}, follows that $\F_{r,s}\in \ucl(\F_r )$.
 \hfill $\blacksquare$
\newline
\medskip
\emph{\textbf{\underline{Proof of Corollary \ref{cor:41}}}}

In the proof of Theorem \ref{thm:II-ucl} above we have already
mentioned that every Lie algebra from $\ucl(\F_r )$ is a $\Phi
_r$-algebra. We, therefore, are to show the converse. Let $B$ be
an arbitrary $\Phi _r$-algebra. Since the formulas of the
collection $\Phi _r$ are universal, we conclude that every
finitely generated subalgebra $A$ of $B$ is a $\Phi _r$-algebra
and consequently lies in $\ucl(\F_r )$. Therefore, $B \in
\ucl(\F_r )$. \hfill $\blacksquare$

The proofs of Theorem \ref{thm:II-F-ucl} and Corollary
\ref{cor:II-F-ucl} are analogous to the ones of Theorem
\ref{thm:II-ucl} and Corollary \ref{cor:II-ucl}.

\begin{rem}
In {\rm \cite{AGFMLA1}} the authors show (see Lemma {\rm 4.8})
that for any torsion free $R$-module $M$ the algebra $\F_{r}\oplus
M$ $\F_{r}$--embeds into the algebra $\F_{r,s}$ for some $s \in
\mathbb{N}$. In what follows that $\Phi' _r$-algebras can be
treated as $\F_{r}$--subalgebras of $\F_{r,s}$.
\end{rem}
\medskip
\emph{\textbf{\underline{Proof of Theorem \ref{thm:43}}}}

The statements of the theorem regard a universal class of a single
object. From general model-theoretical facts follows that it
suffices to prove the first statement only.

All the series of axioms $\Phi 1- \Phi 6$ are obviously recursive,
provided that the field $k$ is finite. We now treat the series
$\Phi 7$. The axioms of this series are enumerated by the set of
finite systems of equations over the  module $\Fit(\F_r)$ which
are inconsistent over $\Fit_{\Delta}(\F_r)$. According to
\cite{Seid}, the compatibility problem for systems of equations
over finitely generated modules over Noetherian commutative rings
is decidable. Since $\Fit_{\Delta}(\F_r)$ is a finitely generated
module over Noetherian commutative ring $R_{\Delta}$ and since the
set of all finite systems of equations over $\Fit(\F_r)$ is
recursive (in the event that the ground field $k$ is finite), the
statement follows. \hfill $\blacksquare$

The proof of Theorem \ref{thm:44} is analogous and, therefore,
omitted.
\medskip
\newline
\emph{\textbf{\underline{Proof of Theorem \ref{thm:45}}}}

In the end of Subsection \ref{sec:22} we have pointed out that to
solve a system of equations over $\F_r$ it is convenient to part
it into two. The first one is a linear system of equations over
$k$. Every solution of this system gives rise to a module system
of equations over $\Fit(\F_r)$. In the event that the field $k$ is
finite these problems are algorithmically soluble. \hfill
$\blacksquare$

\section{Irreducible Algebraic Sets over $F_r $ and Dimension}
\label{sec:5}

In Section \ref{sec:2} we have introduced the category of
algebraic sets over an arbitrary Lie algebra and the category of
coordinate algebras. We show there that these two categories are
equivalent. In particular, the classification of coordinate
algebras gives us a classification of algebraic sets. Furthermore,
the algebraic set is defined by its coordinate algebra up to
isomorphism in the respective category. In Subsection \ref{sec:51}
below we give a classification of irreducible algebraic sets over
$F_r$, using the classification of its coordinate algebras and in
Subsection \ref{sec:52}, following the custom of classical
algebraic geometry, we introduce its counterpart -- the definition
of a dimension of an algebraic set. In Theorem \ref{thm:II-dim} we
show how one can find the dimension of an arbitrary irreducible
algebraic set over $F_r$.

\subsection{Classification of Irreducible Algebraic Sets over
$F_r$} \label{sec:51}

Throughout this subsection we use the notation  $a_1 ,\ldots,a_r $
for the free base of the free metabelian Lie algebra $F_r $, $r
\ge 2$. Theorem \ref{thm:222} states that the collection of all
coordinate algebras of irreducible algebraic sets over $F_r $
coincides with the family of all finitely generated algebras from
$F_r - \ucl(F_r )$. Theorem \ref{thm:42} gives a structural
description of such algebras. The following proposition combines
these results.

\begin{prop}[on irreducible coordinate algebras over $\F_r$] \quad
\quad \label{prop:II-irrcoralg} \label{prop:511} Let $\Gamma$ be a
$\F_r$--Lie algebra. Then $\Gamma$ is the coordinate algebra of an
irreducible algebraic set over $\F_r $ if and only if \/ $\Gamma$
is $\F_r $--isomorphic to $\F_r \oplus M$, where here $M$ is a
torsion free module over $k\left[x_1 ,\ldots,x_r \right]$.
\end{prop}

Using Proposition \ref{prop:511}, we next approach the problem of
classification of irreducible algebraic sets over $F_r $.

We shall make use of the following notation.  Let $R = k\left[x_1
,\ldots,x_r \right]$ and let $M$ be a finitely generated
torsion-free module over $R$. Let $\hom_R (M,\Fit(\F_r ))$ be the
set of all $R$-homomorphisms from  $M$ to $\Fit(\F_r )$ treated as
a module over $R$. This set of homomorphisms can be treated from
another point of view. Fix a system of generators $\left\{ m_1
,\ldots,m_n \right\}$ of a module $M$. Defining an
$R$-homomorphism is equivalent to defining the images of the
elements $ m_1 ,\ldots,m_n$. This defines an embedding
\begin{gather} \notag
\begin{split}
\alpha: \hom_R (M,\Fit(\F_r )) \rightarrow  \underbrace
{\Fit(\F_r ) \times \cdots \times  \Fit(\F_r )}_n, & \\
\alpha (\phi )& = (\phi (m_1 ),\ldots, \phi (m_n )),
\end{split}
\end{gather}
 where here $\phi  \in \hom_R (M,\Fit(\F_r ))$. Every
relation imposed on the $n$-tuple $\left\{ m_1 ,\ldots,m_n
\right\}$ is a relation between $\phi (m_1 ),\ldots,\phi (m_n )$.
We, therefore, identify the set  $\hom_R (M,\Fit(\F_r ))$ with its
image $\alpha (\hom_R (M,\Fit(\F_r )))$ in $\underbrace {\Fit(\F_r
) \times \cdots \times \Fit(\F_r )}_n$.

\begin{lem} \label{lem:512}
In this notation, the following chain of one-to-one correspondences
holds {\rm
$$
\hom_R (M,\Fit(\F_r )) \leftrightarrow \hom_{\F_r } (\F_r \oplus
M,\F_r ) \leftrightarrow Y,
$$}
where here $Y$ is an irreducible algebraic set over $\F_r $ such
that $\Gamma(Y)=\F_r  \oplus M$.
\end{lem}
\proof The result follows directly from Lemma  \ref{lem:211} and
Lemma 4.6 in \cite{AGFMLA1}. \hfill $\blacksquare$

\begin{thm}[on irreducible algebraic sets over $\F_r $] \quad  \label{thm:513}
\\
Every irreducible algebraic set over $\F_r$ is,
up to isomorphism, either
\begin{itemize}
    \item a point or
    \item  {\rm $\hom_R(M,\Fit(\F_r ))$} for some finitely generated torsion free
    module  $M$ over the ring $R$.
\end{itemize}
Conversely, any set of the above form is algebraic.
\end{thm}
\proof Consider an arbitrary irreducible algebraic set $Y$ over
the algebra $F_r $. By Proposition \ref{prop:511} its coordinate
algebra $F_r $ is isomorphic to the algebra of the form $F_r
\oplus M$, where $M$ is a finitely generated torsion-free module
over the ring $R$.

If $M = 0$, i.e.  the coordinate algebra is  $F_r $, in which case
there exists the only $F_r$--homomorphism $\varphi $ from $\Gamma
(Y)$ into $F_r $. Set $\varphi (x_i ) = b_i $, $i = 1, \ldots ,n$.
Then, by Lemma \ref{lem:211}, $Y = \left\{ (b_1 ,\ldots,b_l
)\right\} \cong\left\{ 0\right\}$ is a point. Conversely, every
point is an algebraic set, obviously an irreducible one. The
corresponding coordinate algebra is isomorphic to $F_r$ (see
Examples \ref{expl:1} and \ref{ex:coordalg}).

Otherwise $M \ne 0$. On account of Lemma \ref{lem:512} it suffices
to show that the set $\hom_R (M,\Fit(F_r ))$ is an algebraic set
over $F_r $ and that its coordinate algebra is  $F_r $--isomorphic
to the algebra $F_r  \oplus M$.

Let $M=\left< m_1 ,\ldots,m_n| s_1, \ldots ,s_l \right>$ be a
presentation of $M$. The system of module equations $S=
\left\{s_1=0, \ldots ,s_l=0 \right\}$ gives rise to a system of
equations $S_1$ over $\F_r$ (see \cite{AGFMLA1}). Set
$$
S'=S_1 \cup \left\{a_1 a_2 x_i = 0 ,\ i = 1,\ldots,n\right\}.
$$ By Lemma \ref{lem:32}, $V_{F_r }
(S') \subseteq \Fit^n (F_r )$. Consequently, $V_{F_r } (S') =
\hom_R (M,\Fit(F_r ))$. It follows, therefore, that the set
$\hom_R (M,\Fit(F_r ))$ is algebraic.

We show next that $\Gamma (S') \cong _{F_r } F_r \oplus M$. Set
$\theta $ to be the following $F_r $--homomorphism:
$$
\theta :(F_r )_\M \left[X\right] \rightarrow F_r  \oplus M,\quad
\theta (x_i ) = m_i ,\;i = 1,\ldots,n,\quad \theta (a) = a,\;a \in
F_r,
$$
and show that $\ker \,\theta  = \rad (S')$. Here $(F_r)_\M [X]$ is
the free $F_r $--algebra generated by the alphabet $X$ in the
variety of all metabelian Lie algebras $\M$ and the radical
$\rad(S')$ is the radical in $\M$. As mentioned in Section
\ref{sec:2} the definition of a coordinate algebra $\Gamma (S')$
carries over to the variety $\M$.

Finally, we show that $\ker\,\theta  = \rad(S')$. Take a
polynomial $f \in \ker\,\theta $ and write it in the following
form
$$
 f = c+x_1 h_1+\cdots + x_n h_n  + g(x_1 ,\ldots,x_n ),\quad c \in F_r ,\;
 h_i  \in R,\;g(x_1 ,\ldots,x_n ) \in I_X^2 ,
$$
(see Section \ref{sec:2}). Since $F_r  \oplus M$ is the direct
module extension of $F_r $ by the module $M$, we have $c = 0$.
Next, for
$$
g(m_1 ,\ldots,m_n ) = 0 \hbox{ and }  f(m_1 ,\ldots,m_n ) = 0,
$$
we obtain the following relation in the module $M$:
$$
     m_1 h_1  + \cdots + m_n h_n  = 0.
$$

Now, it is essentially immediate that for any point
$$
(b_1,\ldots,b_n ) \in V_{F_r } (S')
$$
holds:
$$
     b_1 h_1  + \cdots + b_n h_n  = 0,
$$
i.e.
$$
x_1 h_1  + \cdots + x_n h_n  \in \rad(S').
$$
Therefore,  $f \in \rad(S')$.

Suppose next that $f \in \rad(S')$ and show that $f \in
\ker\,\theta $. Since  $(0,\ldots,0) \in V_{F_r } (S')$, we have
$c = 0$. Moreover,
$$
     g(m_1 ,\ldots,m_n ) = 0.
$$
Verify that
$$
     m_1 h_1  + \cdots + m_n h_n  = 0.
$$
Assume the converse, then there exists an  $R$-homomorphism $\phi
\in \hom_R (M,\Fit(F_r ))$ so that
$$
     \phi (m_1 h_1  + \cdots + m_n h_n ) \ne 0
$$
(see Lemma 4.7 in \cite{AGFMLA1}). This derives a contradiction,
for $f \in \rad(S')$. \hfill $\blacksquare$

\begin{rem}
The system $S'$, associated to $\F_r  \oplus M$ depends on the
choice of a presentation of the module $M$. We term such systems
\emph{canonical} for the algebra $\F_r  \oplus M$.
\end{rem}

\begin{cor}
In the one-dimensional affine space $\F_r^n $, $n = 1$ every
irreducible algebraic set is, up to isomorphism, either
\begin{itemize}
    \item a point or
    \item  {\rm $\Fit(\F_r )$}.
\end{itemize}
\end{cor}
\proof The coordinate algebra of an irreducible algebraic set is
$F_r $--isomorphic to the algebra $F_r  \oplus M$. If $M = 0$ then
the algebraic set is a point. If  $M \ne 0$, since $M$ is
one-generated, the module is the free module $T_1 $. In which case
$\hom_R (T_1 ,\Fit(F_r ))$  is isomorphic to $\Fit(F_r )$. \hfill
$\blacksquare$

\subsection{Dimension} \label{sec:52}

\begin{defn}
Let $Y$ be an irreducible algebraic set. As is the custom in
algebraic geometry, a maximum of all integers $m$ such that there
exists a chain of irreducible algebraic sets
$$
Y = Y_0  \varsupsetneq Y_1  \varsupsetneq \ldots \varsupsetneq
Y_m.
$$
is termed the dimension of $Y$ and is denoted by $\dim (Y)$.
\end{defn}

\begin{defn}
Let $Y = Y_1  \cup \ldots \cup Y_l$ be an expression of an
algebraic set $Y$ (not necessarily irreducible) as a finite union
of irreducible algebraic sets (see Theorem \ref{thm:214}). We
define the dimension of $Y$ (and denote it by $\dim (Y)$) to be
the maximum of the dimensions of all irreducible components.
\end{defn}

Let $Y$ be an irreducible algebraic set over  $\F_r $. Let $\Gamma
(Y) \cong _{\F_r } \F_r  \oplus M$ for some finitely generated
torsion free module $M$ over $R = k\left[x_1 ,\ldots,x_r \right]$.
Theorem \ref{thm:II-dim} shows that $\dim (Y)$ is uniquely defined
by the module $M$. Recall that the rank $r(M)$ of the module $M$
over the ring $R$ is a maximum of cardinalities of linearly
independent over $R$ sets of elements from $M$.
 By the definition we set $r(\Gamma (Y))=r(M)$.

\begin{thm} \label{thm:II-dim}
For an  irreducible algebraic set $Y$ over $\F_r$ holds
$$
\dim (Y) = r(\Gamma (Y)) = r(M).
$$
\end{thm}
\proof Let
$$
Y_m  \varsubsetneq \ldots\varsubsetneq Y_1 \varsubsetneq Y_0  = Y
$$
be a strictly descending chain of irreducible algebraic sets.
Obviously $Y_m$ is a point, $\dim (Y_m) = 0$ , $ r(\Gamma (Y_m ))
= 0$. By Lemma \ref{lem:213} the inclusion $Y_{i + 1}
\varsubsetneq Y_i $ induces an $\F_r$--epimorphism of respective
coordinate algebras $\varphi: \F_r \oplus M_i \rightarrow \F_r
\oplus M_{i + 1}$, moreover $\ker \varphi \ne 0$. According to
Lemma 4.8 in \cite{AGFMLA1}, $\varphi$ induces an $R$-epimorphism
$\phi :M_i \rightarrow M_{i + 1} $, where $\ker \phi  \ne 0$.
Consequently, $r(M_i ) > r(M_{i + 1})$. In what follows that $\dim
(Y) \le r(\Gamma (Y))$.

To prove the converse inequality suppose that $r(\Gamma (Y)) =
r(M) = n$. Let $N$ be the isolated submodule generated by a
nonzero element $0 \ne m \in M$, and let $\phi$  be the canonical
$R$-epimorphism from $M$ onto $M_1 = {\raise0.7ex\hbox{$M$}
\!\mathord{\left/
 {\vphantom {M N}}\right.\kern-\nulldelimiterspace}
\!\lower0.7ex\hbox{$N$}}$. Then $M_1 $ is a torsion free module
over $R$ and $r(M_1) = n - 1$. Applying the reverse of the above
argument we conclude that $\dim (Y) \ge n$. \hfill $\blacksquare$

\newpage
\textit{Daniyarova Evelina Yur'evna},
\newline
644043, Russia, Omsk, Spartakovskaya st. 13-8,
\newline
tel. +7 3812 232239,
\newline
\textsl{ e-mail:} \verb"evelina_om@mail333.com"
\newline
Omsk Branch of Institute of Mathematics
\newline
(Siberian branch of Russian Academy of Science)

\bigskip

\textit{Kazatchkov Ilia Vladimirovich},
\newline
644046, Russia, Omsk, Pushkin st. 136-22,
\newline
tel. +7 3812 312315,
\newline
\textsl{e-mail:} \verb"kazatchkov@mail333.com"
\newline
Omsk Branch of Institute of
Mathematics
\newline
(Siberian branch of Russian Academy of Science)

\bigskip

\textit{Remeslennikov Vladimir Nikanorovich},
\newline
644099, Russia, Omsk, Ordjonikidze st. 13-202,
\newline
tel. +7 3812 240914,
\newline
\textsl{e-mail:} \verb"remesl@iitam.omsk.net.ru"
\newline
Omsk Branch of Institute of Mathematics
\newline
(Siberian branch of Russian Academy of Science)

\end{document}